\documentclass[10pt]{amsart}
\usepackage{amsfonts,amssymb,verbatim}
\usepackage[scriptlabels,height=8mm,width=8mm]{diagrams}
\diagramstyle[PostScript=dvips]
\input xy
\xyoption{all}

\newarrow{Equal}{}{=}{}{=}{}
\newcommand\dis{\displaystyle}

\newcommand{\cale}{{\mathcal E}}
\newcommand{\calf}{{\mathcal F}}

\newcommand{\calh}{{\mathcal H}}
\newcommand{\cali}{{\mathcal I}}

\newcommand{\caln}{{\mathcal N}}
\newcommand{\calo}{{\mathcal O}}
\newcommand{\calp}{{\mathcal P}}

\newcommand{\calz}{{\mathcal Z}}

\newcommand{\calox}{{\calo_X}}

\newcommand{\caloz}{{\calo_Z}}

\newcommand{\bbbc}{{\mathbb C}}

\newcommand{\bbbh}{{\mathbb H}}
\newcommand{\bbbl}{{\mathbb L}}

\DeclareMathOperator{\At}{At}
\DeclareMathOperator{\at}{at}

\DeclareMathOperator{\ch}{ch}

\DeclareMathOperator{\Ext}{Ext}
\DeclareMathOperator{\Hom}{Hom}
\DeclareMathOperator{\Hilb}{Hilb}
\DeclareMathOperator{\HH}{H}
\DeclareMathOperator{\id}{id}

\DeclareMathOperator{\ob}{ob}

\DeclareMathOperator{\tr}{tr}

\newcommand{\lr}{\,\left\lrcorner\right.\,}

\newcommand{\bdot}{\mathop{\mbox{\boldmath $\cdot$}}}

\newcommand{\TX}{\Theta_X}

\theoremstyle{definition}
\newtheorem{definition}{Definition}[section]
\newtheorem{sit}[definition]{}

\theoremstyle{plain}
\newtheorem{proposition}[definition]{Proposition}
\newtheorem{theorem}[definition]{Theorem}

\begin{document}

\title[Atiyah-Chern Character and Semiregularity]%
{The Atiyah-Chern Character yields the Semiregularity Map
as well as the Infinitesimal Abel-Jacobi Map}

\author[R.-O.~Buchweitz]{Ragnar-Olaf Buchweitz}
\address{Dept. of Mathematics, University of Toronto,
Toronto, Ont. M5S 3G3, Canada}
\email{ragnar@math.utoronto.ca}
\author[H.~Flenner]{Hubert Flenner}
\address{Fakult\"at f\"ur Mathematik, Ruhr-Universit\"at Bochum,
Geb.~NA 2/72, 44781 Bochum, Germany}
\email{Hubert.Flenner@ruhr-uni-bochum.de}

%
\thanks{R.O.B.~was partly supported by NSERC grant 3-642-114-80.
\endgraf 
Both authors were partly supported by a grant from the Volkswagen Foundation 
under the Research in Pairs program at Math.~Forschungsinstitut Oberwolfach.}
\date{\today}
\maketitle

\section*{Introduction} 

The purpose of this work is to construct a general semiregularity map 
for cycles on a complex analytic or algebraic manifold and to show 
that such semiregularity map can be obtained from the classical tool 
of the Atiyah-Chern character.  The first part of the paper is fairly 
detailed, deducing the existence and explicit form of a generalized 
semiregularity map from known results and constructions.  As a 
corollary we obtain in the second part as well a description of the 
infinitesimal Abel-Jacobi map for smooth cycles as the leading term of 
this generalized semiregularity map, indicate why for locally complete 
intersections the appropriate component of our semiregularity map 
coincides with the one constructed by Bloch \cite{Bloch}, and give 
applications to embedded deformations and deformations of coherent 
modules.

At the request of several participants at this conference, we restrict 
ourselves mainly to the case of a smooth ambient space, avoiding thus 
the formidable technical machinery of resolvents, powers of cotangent 
complexes, traces and so on.  The underlying ``metaresult'' in its 
utmost generality is only stated at the end with details to be given 
in \cite{BF}.

We wish finally to thank the editorial committee of these proceedings for 
their inexhaustible patience awaiting completion of this paper.

\section{Chern Character and Griffiths' Transversality}

\begin{sit}
\label{formula}
	Let $Z \subseteq X$ be a codimension $q$ cycle on a compact 
	K\"ahler manifold or a smooth quasi-projective variety $X$ defined 
	over a field of characteristic zero.  With $\Omega^p_X$ the sheaf 
	of (holomorphic or algebraic) differential forms of order $p$ on 
	$X$, the (cohomological) cycle or fundamental class of $Z$ in 
	$\HH^{q}(X, \Omega^q_X)$ has the following classical description 
	due to Grothendieck \cite[p.151]{Gr1}, see also \cite[2.18]{Murre} 
	for a recent treatment.
	
	Let $c_{p}(\calo_{Z}), \ch_{p}(\calo_{Z})\in \HH^{p}(X, \Omega^p_X)$ 
	denote the $p^{\rm th}$ {\em Chern class\/} of $Z$ and 
	the $p^{\rm th}$ component of the {\em Chern character\/} of $Z$, respectively.
	
	The {\em fundamental (Hodge) class of $Z$\/} is
	$$
	[Z] = \frac{(-1)^{q-1} c_q(\calo_Z)}{(q-1)!} = \ch_q (\caloz) \in 
	\HH^q(X, \Omega^q_X)\,.
	$$
	As the components $\ch_{p}(\calo_{Z})$ vanish for $p$ below the 
	codimension, the fundamental class is thus the leading term of 
	the Chern character, $\ch(\caloz) =\sum_{p\ge q} 
	\ch_{p}(\caloz)\in \prod_{p} \HH^{p}(X, \Omega^{p}_X)$.
\end{sit}

\begin{sit}
	Let $X$ be an analytic space or a scheme over a field.  If 
	$X$ is {\it smooth\/}, its first order (flat) deformations
	are parametrized by $\HH^1(X, \TX)$, the first cohomology group 
	of the locally free sheaf of (holomorphic or algebraic) vector 
	fields $\TX= \calh om_{X}(\Omega^1_X,\calox)$.
	
	The group $\HH^1(X,\TX)\cong 
	\Ext^{1}_{X}\left(\Omega^1_X,\calox\right)$ acts through 
	graded derivations with respect to the total degree $p+q$ on the algebra 
	$\prod_{p,q}\HH^q(X, \Omega^p_X)$, sending a class of 
	type $(p,q)$ to one of type $(p-1,q+1)$. Denoting the 
	action of the class of a first order deformation $*$ on $c\in \HH^q 
	\left(X, \Omega^p_X\right)$ by $ *\lr c$, as in \cite[III.164, 
	Prop.  10]{ALG} or \cite[6.1]{Bloch}, a codimension $q$ cycle $Z$ 
	on $X$ defines thus a commutative diagram
	$$
	\xymatrix{
	\HH^1(X, \TX)
	\ar[rr]^{* \lr\ch(\caloz)}
	\ar[dr]_{* \lr [Z]}
	&& \prod\limits_{p}\HH^{p+1}(X, \Omega_X^{p-1})
	\ar[dl]^{\rm projection}\\
	& \HH^{q+1}(X, \Omega_X^{q-1})&.}
	$$
The geometric meaning of the action of first order deformations on 
cycle classes is explained by Bloch's interpretation of Griffiths'
Transversality Theorem: 
\end{sit}

\begin{theorem}{\em (\cite[(4.2)]{Bloch}; see \cite[2.3]{Voi} for the 
analytic case.)}\\
	Given $\xi \in \HH^1(X, \TX)$, the unique horizontal lift of $[Z]$ 
	relative to the Gau\ss-Manin connection defined by the 
	infinitesimal deformation $\xi$ stays of type $(q, q)$ if and 
	only if $ \xi \lr [Z] = 0$.
	\qed
\end{theorem}

One way of stating Grothendieck's (infinitesimal) {\em Variational 
Hodge Conjecture\/} \cite{Gr2}, see also \cite[\S 3]{Ste}, is to ask 
whether the obstruction expressed by Griffiths' Transversality is 
algebraic: If the unique horizontal lift of $[Z]$ along a first order 
deformation $\xi$ remains of Hodge type $(q,q)$, is there then some 
cycle representing the cycle class of $[Z]$ that deforms along $\xi$ 
in $X$?  To begin with one may ask, as Bloch did in \cite{Bloch}, 
what is the relationship between the obstruction to deform the given 
$Z$ along $\xi$ and the cohomological obstruction $\xi\lr [Z]\,$?

This leads to the notion of {\em semiregularity\/} that we now recall.

\section{The Structure of the Hilbert Scheme and Semiregularity}

Let $X$ be an arbitrary analytic space or scheme over some field and 
let $Z \subseteq X$ be any closed subspace defined by a sheaf of 
ideals $\cali \subseteq \calox$.  The embedded deformations of 
$Z\subseteq X$ are parametrized in the analytic case by the {\em 
Douady space} $H_{X}$ at the point corresponding to $Z\subseteq X$, 
and in the algebraic case by the {\em Hilbert scheme\/} $\Hilb_{X}$ at 
$Z$.

\begin{sit}
	The tangent space to $H_X$ or $\Hilb_X$ at the point corresponding to
	$Z$ is
	$$
	T^{1}_{Z/X} = \Hom_Z(\cali/\cali^2, \caloz) \cong
	\HH^0(Z, \caln_{Z/X})
	$$
	where $\caln_{Z/X}$ is the normal sheaf to $Z$ in $X$.
	The obstruction space is contained in $T^2_{Z/X}$, the second 
	cohomology group of the relative cotangent complex of $Z$ in $X$,
	see e.g. \cite{Ill},\cite{Rim},\cite{Fle}, \cite{Pal}. One 
	has the following inclusions and equalities that approximate this group
	and describe it in many situations:
	$$
	\begin{array}{lllll}
	\HH^1(Z, \caln_{Z/X}) & \subseteq 
	      & \Ext^1_Z(\cali/\cali^2, \caloz)
	      & \subseteq  &T^2_{Z/X}\\
	& \uparrow  & & \uparrow & \\
	&=& \text{if $Z$ is locally} &=&\text{if $Z$ is generically a complete} \\
	&&\text{a complete intersection}&&\text{intersection and has no}\\
	&&&&\text{embedded 
	components.}
	\end{array}
	$$
\end{sit}

\begin{sit}
	With $X$ still an arbitrary space, denote $T^{1}_{X}$ its space of 
	first order deformations.  There is then a natural obstruction map 
	$\ob_{Z/X}\colon T^{1}_{X}\to T^2_{Z/X}$ so that $Z$ can be 
	deformed along a given (first order) deformation $\xi$ of $X$ iff 
	$\ob_{Z/X}(\xi)=0$.
	
	If $X$ is smooth, then $T^{1}_{X}\cong \HH^{1}(X,\TX)$ and the 
	obstruction map factors through the restriction map
	$\HH^{1}(X,\TX) \xrightarrow{\quad}
	\HH^{1}(Z,\TX\otimes\caloz)\cong \HH^{1}(Z,\calh 
	om_{Z}(\Omega^{1}_X\otimes\caloz,\caloz))$ followed by 
	$\HH^{1}(Z,\ \,)$ of the $\caloz$-dual of the 
	Jacobi map $j:\cali/\cali^2\to \Omega^{1}_X\otimes\caloz$,
	$$
	\ob_{Z/X}\colon T^{1}_{X}\cong \HH^{1}(X,\TX) \xrightarrow{\quad}
	\HH^{1}(Z,\TX\otimes\caloz)\xrightarrow{\HH^{1}(Z,j^{*})}
	\HH^1(Z, \caln_{Z/X}) \subseteq 
	T^2_{Z/X}\,.
	$$
\end{sit}
	
\begin{sit}	
	The question referred to above is now: Is there a natural {\em 
	semiregularity map\/}
	$$
	\sigma: T^2_{Z/X} \xrightarrow{\quad} \HH^{q+1}(X, \Omega^{q-1}_X)
	$$
	such that the diagram
	$$
	\xymatrix{
	\HH^1(X, \TX)
	\ar[rr]^{\ob_{Z/X}(*)}
	\ar[dr]_{*\lr [Z]}
	&& T^2_{Z/X}
	\ar@{.>}[dl]^\sigma\\
	& \HH^{q+1}(X, \Omega^{q-1}_X)}
	$$
	commutes?
\end{sit}

\begin{sit}
	If a semiregularity map exists, and if it is furthermore {\em 
	injective}, then $Z$ is called {\em semiregular} in $X$.  
	
	Such a semiregularity map was first introduced by Severi 
	\cite{Sev} for a curve on a surface and he coined the term 
	{semiregular\/} in that context.  Kodaira-Spencer \cite{K-S} later 
	generalized Severi's construction and results to arbitrary 
	divisors on a complex manifold using the tools from the 
	deformation theory they had just developed.  Bloch \cite{Bloch} 
	extended the notion to arbitrary cycles on smooth complex 
	projective schemes, constructed a semiregularity map for cycles 
	that are locally complete intersection, showed that the Hilbert 
	scheme $\Hilb_{X}$ is smooth at the point corresponding to a 
	semiregular locally complete intersection $Z\subseteq X$, and 
	deduced then the variational Hodge conjecture for the 
	corresponding cycle class $[Z]$.
	
	We will show here that a semiregularity map always exists, and 
	that its construction could have been achieved more than 25 years ago, by 
	anyone reading \cite{Ill} and \cite{Bloch} side by side!
	
	The key ingredient is to forget part of the 
	structure and to construct a more general semiregularity map for 
	coherent sheaves on $X$.  The underlying idea is simple: The 
	fundamental class is not only defined for a cycle in $X$ but for 
	any element of ${K^{0}}(X)$, the Grothendieck group of the smooth 
	space $X$.  Thus one may ask for a semiregularity map for any 
	representative of a class in that group, that is for any coherent 
	sheaf on $X$, or, if one allows $X$ to become singular, for any 
	{\em perfect complex} of sheaves on $X$.
\end{sit}

\section{Deformations of Modules and the Atiyah-Chern Character}

If $Z\subseteq X$ is a closed subspace of $X$, we may consider instead 
of embedded deformations of $Z$ the deformations of $\calo_{Z}$ as 
$\calo_{X}$-module.  This is a forgetful functor: a deformation of 
$\calo_{Z}$ as module will only be an invertible sheaf on its support, not 
necessarily its structure sheaf.  

\begin{sit}
\label{forget}
	For any coherent sheaf $\calf$ on 
	$X$, the space of first order deformations is given by $\Ext^1_X (\calf, 
	\calf)$ with obstructions in $\Ext^2_X (\calf, \calf)$, and so there are
	accordingly natural forgetful transformations
	\begin{equation}
		f^{i}\colon T^i_{Z/X} \xrightarrow{\quad}\Ext^i_X (\caloz, \caloz)\,.
		\tag{\ref{forget}.1}
	\end{equation}	
	These maps can be explained quite simply within the theory of 
	resolvents and cotangent complexes, as developed in \cite{Fle}, 
	\cite{Ill}, \cite{Pal} based on  \cite{Qui}: a 
	resolvent for $Z$ over $X$ is a differential graded free $\calox$-algebra that 
	resolves $\caloz$, the groups $T^i_{Z/X}$ constitute the 
	hypercohomology groups of the $\calox$-derivations of said 
	resolvent, whereas the groups $\Ext^i_X (\caloz, 
	\caloz)$ represent the hypercohomology of its $\calox$-linear 
	endomorphisms.  Now any $\calox$-derivation is $\calox$-linear by 
	definition, whence the comparison maps; see \cite{BF} for further details.
\end{sit}

\begin{sit}
\label{Z2}	
	Note that $\Ext^i_X (\caloz, \caloz)$ always contains 
	$\HH^{i}(X,\caloz) = \HH^{i}(Z,\caloz)$ naturally as a direct summand and that
	$f^{i}$ always maps into the complement. For $i=1$ one has indeed 
	a split exact sequence
	$$
	0\to  T^1_{Z/X}=\HH^{0}(X,\caln_{Z/X})\xrightarrow{f^{1}}\Ext^1_X (\caloz, \caloz)
	\to \HH^{1}(Z,\caloz)\to 0
	$$ 
	and $\HH^{1}(Z,\caloz)$ is the tangent space to the Picard scheme 
	of $Z$ at $\caloz$. This formalizes that a first order deformation 
	of $\caloz$ as $\calox$-module splits into an embedded deformation of 
	$Z$ and a deformation of the invertible sheaf $\caloz$ on $Z$. 

	When restricted to 
	$\Ext^{i-1}_{Z}(\cali/\cali^{2},\caloz)\subseteq T^i_{Z/X}$, the 
	forgetful maps $f^{i}$ become particularly explicit: The extension 
	defined by the first fundamental neighbourhood of $Z$ in $X$,
	\begin{equation*}
	\calz^{(2)}\quad\equiv\quad 0\xrightarrow{\quad} \cali/\cali^{2} 
	\xrightarrow{\quad} \calox/\cali^{2}\xrightarrow{\quad}\caloz 
	\xrightarrow{\quad} 0\,,
	\end{equation*}
	defines a connecting homomorphism 
	$\Ext^{i-1}_{X}(\cali/\cali^{2},\caloz)\to \Ext^i_X (\caloz, 
	\caloz)$ which composed with the forgetful map 
	$\Ext^{i-1}_{Z}(\cali/\cali^{2},\caloz)\to 
	\Ext^{i-1}_{X}(\cali/\cali^{2},\caloz)$ yields the restriction of 
	$f^{i}$.
\end{sit}

\begin{sit}
	Returning to the context of coherent $\calox$-modules on a smooth space $X$ we 
	want thus to construct a general semiregularity map for coherent 
	$\calox$-modules $\calf$,
	$$
	\tau^{2}=\left(\tau^{2,p}\right)_{p\ge 0}\colon\Ext^2_X (\calf, 
	\calf)\xrightarrow{\quad} \prod_{p\ge 0} \HH^{p+2} \left(X, 
	\Omega^p_X\right)\,,
	$$
	such that the desired semiregularity map for a cycle $Z$ of 
	codimension $q$ will be
	$$
	\sigma= \tau^{2,q-1} f^{2}\colon 
	 T^2_{Z/X} \xrightarrow{\quad}\Ext^2_X (\caloz, \caloz)
	  \xrightarrow{\quad}  \HH^{q+1}(X, \Omega^{q-1}_X)\,.
	$$
	A natural such map $\tau$ has been known for a long time! 
\end{sit} 

\begin{sit}
	To describe it, recall the definition of the {\em Atiyah class\/} 
	of an $\calox$-module $\calf$: Let $X^{(2)}\subseteq X\times X$ 
	denote the first fundamental neighbourhood of the diagonal 
	$\Delta:X\hookrightarrow X\times X$ and denote $p_{i}\colon 
	X^{(2)}\to X$ the projection onto the $i^{\rm th}$ factor.  The 
	sheaf of $1$-jets or principal parts of first order associated to 
	$\calf$ is defined to be $\calp^{1}(F) = p_{2*}p^{*}_{1}\calf$.  
	Restricting $\calp^{1}(F)$ to the diagonal itself defines a 
	canonical exact sequence of (right) $\calox$-modules, the {\em 
	Atiyah sequence\/} of $\calf$,
	$$
	\At_{X}(\calf)\quad\equiv\quad 0\xrightarrow{\quad} \calf\otimes 
	\Omega^1_X\xrightarrow{\quad} \calp^{1}(\calf) \xrightarrow{\quad}\calf 
	\xrightarrow{\quad} 0
	$$
	whose class
	$$
	\at_{X}(\calf)\in \Ext^{1}_X (\calf, \calf \otimes \Omega^1_X)
	$$
	is the (first) {\em Atiyah class\/} of $\calf$. 
\end{sit}

\begin{sit}
    Combining the composition or Yoneda product with the exterior product 
	on $\Omega^{\bdot}_{X}$ defines on the vector space 
	$\prod_{i,j}\Ext^{i}_{X}(\calf,\calf\otimes \Omega^j_X)$ an algebra 
	structure with multiplication
	$$
	\smallsmile\colon \Ext^{a}_{X}(\calf,\calf\otimes \Omega^b_X)\times
	\Ext^{c}_{X}(\calf,\calf\otimes \Omega^d_X)\xrightarrow{\quad} 
	\Ext^{a+c}_{X}(\calf,\calf\otimes \Omega^{b+d}_X)\,.
	$$
	One can thus form powers of the Atiyah class,
	$$
	\at^{p}_{X}(\calf)= \at_{X}(\calf)^{\smallsmile p}
	\in \Ext^{p}_X (\calf, \calf \otimes \Omega^p_X)\,,
	$$
	and then, in characteristic zero, the {\em Atiyah-Chern 
	character\/}
	$$
	 e^{\dis \at_{X}(\calf)}=
	\sum_{p\ge 0}\frac{\at^{p}_{X}(\calf)}{p!}
	\in \prod_{p} \Ext^{p}_X (\calf,\calf \otimes \Omega^p_X)\,.
	$$	
	The (powers of the) Atiyah class and the Atiyah-Chern character 
	are functorial in the appropriate sense; in particular these 
	notions localize to yield local sections in $\prod_{p} \cale 
	xt^{p}_X (\calf,\calf \otimes \Omega^p_X)$.
	Various explicit representatives of these extension classes are 
	worked out in \cite{ALJ} in the algebraic case.
\end{sit}

\begin{sit}
	As $X$ is supposed to be smooth, any coherent $\calox$-module 
	$\calf$ admits locally a finite resolution by locally free 
	sheaves, that is, it represents a {\em perfect complex\/}, and 
	thus there exist natural trace maps
	$$
	\tr^{p}\colon\Ext^{p}_X (\calf,\calf \otimes \Omega^p_X)\to 
	\HH^p\left( X, \Omega^p_X\right) 
	$$ 	
	that uniquely extend the usual trace $\tr^{0}$ of an endomorphism 
	of $\calf$.  The quite intricate technicalities in the general 
	setting were overcome in the algebraic case by Illusie, first in 
	\cite[Exp.~I]{SGA 6}, then in \cite{Ill}, and the analytic case 
	was dealt with explicitly by O'Brian, Toledo, Tong in \cite{OTT}.
\end{sit}

It is a classical result that the trace 
of (the inverse of) the Atiyah-Chern character yields the Chern character:

\begin{theorem}
	{\em (Atiyah \cite{At}; Illusie \cite{Ill}; O'Brian-Toledo-Tong 
	\cite{OTT})\/}\quad 
	If $\calf$ is a coherent sheaf on the smooth space $X$, then
	$$
	\ch(\calf) = \tr\left(e^{\dis-\at_X (\calf)}\right) \in \prod_p
	\HH^p\left( X, \Omega^p_X\right)\,.
	$$
	In particular, for a closed subspace $Z\subseteq X$ of codimension 
	$q$, one has
	$$
	[Z] = \ch_{q}(\caloz) = \frac{(-1)^{q}\tr\left({\at^{q}_X 
	(\caloz)}\right)}{{q!}}\,.
	$$
	\qed
\end{theorem}

The minus sign appears here as traditionally the first Chern class
is the trace of the opposite of the Atiyah class, see \cite[Prop.12]{At}, 
\cite[V.5.4.1, 5.9.4]{Ill}. Note also that extension to the algebraic case 
removed the factor $2\pi i$ originally found in the exponent.

The link between this theorem and the existence of a semiregularity 
map is the following result, due to Illusie \cite[IV.3.1.8]{Ill} in 
much vaster generality.

\begin{proposition}
\label{atiyah is obstruction}
	If $\calf$ is a coherent sheaf on the smooth space $X$, then the 
	obstruction $\ob_{\calf}\colon T^{1}_{X}\to 
	\Ext^{2}_{X}(\calf,\calf)$ to deform $\calf$ (flatly) along a given first order 
	deformation of $X$ can be realized as
	$$
	\ob_{\calf}\colon  T^{1}_{X}\cong \Ext^{1}_{X}(\Omega^1_X,\calox)
	\xrightarrow{\calf\otimes(\ \,)}  \Ext^{1}_{X}(\calf\otimes \Omega^1_X,\calf)
	\xrightarrow{(\ \,){\scriptscriptstyle \circ}(-\dis\at_X (\calf))} 
	\Ext^{2}_{X}(\calf,\calf)\,.
	$$
	\qed
\end{proposition}
In comparison to Illusie's original treatment we changed here the sign in 
front of the Atiyah class to incur fewer signs later on. 

\begin{sit}
	The obstructions for embedded deformations of a subspace $Z$ 
	versus those for deformations of the module $\caloz$ are related 
	explicitly by means of the following commutative diagram whose 
	rows are exact sequences of $\calox$-modules, namely the Atiyah 
	sequence of the $\calox$-module $\caloz$ and the exact sequence of 
	the first infinitesimal neighbourhood of $Z$ in $X$, respectively:
	\begin{diagram}
		\At_{X}(\caloz)&\quad\equiv &0&\rTo&\caloz\otimes 
		\Omega^1_X&\rTo &\calp^{1}(\caloz) &\rTo &\caloz &\rTo & 0\\
		\uTo &&&&\uTo^{j}&&\uTo^{j'}&&\uEqual\\	
		\calz^{(2)}&\quad\equiv &0&\rTo&\cali/\cali^{2}&\rTo& 
		\calox/\cali^{2}&\rTo&\caloz &\rTo& 0\,.
	\end{diagram}
	The map $j'$ is given by $j'(x) = 1\otimes p_{2}^{*}(x)\in 
	\calp^{1}(\caloz) = p_{1}^{*}\caloz\otimes(\calo_{X\times 
	X}/\cali^{2}_{\Delta})$, where $x$ is a local section of 
	$\calox/\cali^{2}$ and $\cali_{\Delta}$ denotes the ideal of the 
	diagonal in $X\times X$.  Indeed, for a local section $f\in \cali$ 
	one has that $j(f)$ lifts to
	\begin{align*}
		{\tilde j}(f) &= 1\otimes df \in \calox\otimes \Omega^1_X\\ 
		     &\mapsto \left (p_{2}^{*}(f) - p_{1}^{*}(f)\right) = 
			     1\otimes f- f\otimes 1 \in \calo_{X\times X}/\cali^{2}_{\Delta}\\
			 &\mapsto 1\otimes p_{2}^{*}(f) \in p_{1}^{*}\caloz\cong 
			     \calo_{Z\times X}/\calo_{Z\times X}\cali^{2}_{\Delta}\\
			 &\mapsto j'(f) \in \calp^{1}(\caloz)\,.
	\end{align*}
	In view of (\ref{Z2}) this means that indeed 
	$\ob_{\caloz}=f^{2}{\scriptscriptstyle \circ}\ob_{Z/X}$ relates the 
	obstruction maps of the two deformation problems.
\end{sit}

Now we are ready to define

\section{The Generalized Semiregularity Map}

\begin{sit}
	The homomorphism of $\calox$-algebras $\calox\to 
	\Omega^{\bdot}_{X}$ induces an algebra homomorphism
	$$\prod_{i}\Ext^{i}_{X}(\calf,\calf)\xrightarrow{\quad}
	\prod_{i,j}\Ext^{i}_{X}(\calf,\calf\otimes \Omega^j_X)$$
	which extends the natural right action of the Yoneda Ext-algebra 
	$\Ext_{X}(\calf,\calf)$ on the target.
	Contraction combined with Yoneda product defines further a left action
	$$
	\Ext^{r}_{X}(\Omega^s_X,\calox)\times 
	\Ext^{i}_{X}(\calf,\calf\otimes \Omega^j_X)\xrightarrow{\ \lr\ }
	\Ext^{r+i}_{X}(\calf,\calf\otimes \Omega^{j-s}_X)  
	$$
	and its restriction to $\Ext^{1}_{X}(\Omega^1_X,\calox)$ acts through 
	graded algebra derivations with respect to total degree,
	$$
	\xi\lr \left(\omega\smallsmile \omega'\right) = \left(\xi\lr \omega\right)\smallsmile 
	\omega' + (-1)^{i+j}\omega\smallsmile\left(\xi\lr \omega'\right)\,,
	$$
	for $\xi\in\Ext^{1}_{X}(\Omega^1_X,\calox), 
	\omega\in \Ext^{i}_{X}(\calf,\calf\otimes \Omega^j_X)$.
\end{sit}

\begin{sit}
	Combining these actions, one obtains for any coherent sheaf 
	$\calf$ a family of {\em generalized semiregularity 
	maps\/},
	$$
	\tau^{i} = \left(\tau^{i,p}\right)_{p\ge 0} =
	\tr\left(e^{\dis-\at_X (\calf)}\smallsmile *\right) \colon
	\Ext^{i}_{X}(\calf,\calf)\xrightarrow{\quad}\prod_{p}
	\HH^{i+p}(X,\Omega^{p}_X) 
	$$
	that is natural in $\calf$ and whose value on $\id_{\calf}\in 
	\Ext^{0}_{X}(\calf,\calf)$ is the Chern character of $\calf$.
\end{sit}

\begin{sit}
	If $Z\subseteq X$ contains the support of $\calf$, then the trace 
	map factors through the local cohomology with support in $Z$, see 
	\cite[V.6]{Ill} or \cite[4.6]{OTT}, whence one may define Chern 
	character and generalized semiregularity maps with support,
	\begin{align*}
		\ch_{Z}(\calf) &= \tr_{Z}\left(e^{\dis-\at_X (\calf)}\right) \in 
		\prod_p \HH_{Z}^p\left( X, \Omega^p_X\right)\,,\\
		\tau_{Z}&= \tr_{Z}\left(e^{\dis-\at_X 
		(\calf)}\smallsmile *\right) \colon 
		\prod_{i}\Ext^{i}_{X}(\calf,\calf)\xrightarrow{\quad}\prod_{i,p} 
		\HH_{Z}^{i+p}(X,\Omega^{p}_X)\,.
	\end{align*}
The key result is now the following.
\end{sit}

\begin{theorem}
\label{key}
	A coherent $\calox$-module $\calf$ on a smooth space $X$ whose 
	support is contained in $Z\subseteq X$ gives rise to a commutative 
	diagram
	$$
	\xymatrix{\HH^1(X,\TX)
	\ar[rr]^{{\rm ob}_{\calf}}
	\ar[dd]_{*\lr \ch(\calf)}
	&& \Ext^2_X (\calf, \calf)
	\ar[dd]^{\tau_{Z}^{2}}
	\ar[ddll]_{\tau^{2}}\\
	\\
	\prod_{p\ge 1}\HH^{p+1}(X, \Omega^{p-1}_X)
	&&
	\ar[ll]^{\rm canonical}
	\prod_{p\ge 1} \HH^{p+1}_{Z}(X, \Omega^{p-1}_X)\,.}
	$$
\end{theorem}

\begin{proof}
	Illusie's result (\ref{atiyah is obstruction}) shows that 
	$\ob_{\calf}(\xi) = \xi\lr -\at_X(\calf)$ 
	for any first order deformation $\xi\in \HH^{1}(X,\TX)$. The 
	class $\at_X(\calf)$ is of total degree two, and so
\begin{align*}
	\xi\lr \ch(\calf) &= 
	\xi\lr \tr\left(e^{\dis-\at_X (\calf)}\right)\\
	& = \tr\left(\xi\lr 
	e^{\dis-\at_X(\calf)}\right)\qquad\text{by naturality of the trace,}\\
	&= \tr\left( e^{\dis-\at_X(\calf)} \smallsmile\left(\xi\lr 
	-\at_X(\calf)\right) \right) \\
	& \qquad\text{as $\tr(ab)=\tr(ba)$ for even classes $a$ and $b$, 
	whence}\\
	&\qquad\text{$\tr(D(e^{-y})) = \tr(e^{-y}D(-y))$ for any (graded) 
	}\\
	&\qquad\text{derivation $D$ and any (even) element $y$,}\\
	& = \tr\left( e^{\dis-\at_X(\calf)} \smallsmile \ob_{\calf}(\xi)\right)
	\qquad\text{by (\ref{atiyah is obstruction}),}\\
	&=\tau^{2}\left(\ob_{\calf}(\xi)\right)\,.
\end{align*}
\end{proof}

Factoring the obstruction map $\ob_{\caloz}$ through $f^{2}$ gives now the desired 
generalized semiregularity map for cycles.

\begin{theorem} If $Z \subseteq X$ is a closed (analytic or 
algebraic) subspace of the smooth space $X$, then the diagram
$$
\xymatrix{\HH^1(X, \TX)
\ar[rrrr]^{{\rm ob}_{\caloz}}
\ar[drr]^{{\rm ob}_{Z/X}}
\ar[ddrr]_{\lr\negmedspace\ch_{Z}(\caloz)}
\ar@/_2pc/[dddrr]_{\lr\negmedspace\ch(\caloz)}
&&&& \Ext^2_X (\caloz, \caloz)
\ar[ddll]^{\tau^{2}_{Z}}
\ar@/^2pc/[dddll]^{\tau^{2}}
\\
&& T^2_{Z/X} 
\ar[rru]^{f^{2}} 
\ar[d]^{\sigma_{Z}^{2}}&\\
&&\dis\prod_{p\ge 1} \HH_Z^{p+1}(X, \Omega_X^{p-1})
\ar[d]^{\rm can.}\\
&& \dis\prod_{p\ge 1} \HH^{p+1}(X, \Omega_X^{p-1})}
$$
commutes, where we define $\sigma_{Z}^{2}= \tau^{2}_{Z}{\cdot} f^{2}$.

If $Z$ is of codimension $q$, then the component $\sigma_{Z}^{2,q-1}= \tau^{2,q-1}_{Z}{\cdot} f^{2}$ is a semiregularity map with support for 
$Z$, whereas its composition with the canonical map yields the desired 
(absolute) semiregularity map
$$
\sigma^{2,q-1} = \tau^{2}{\cdot} f^{2}\colon 
T^2_{Z/X}\xrightarrow{\quad} \HH^{q+1}(X, 
\Omega_X^{q-1})\,.
$$
\qed
\end{theorem}

\begin{sit}
	More generally, we can define as well a family of generalized 
	semiregularity maps for $Z$ (with support) through
	$$
	\sigma_{(Z)}^{i} = \tau^{i}_{(Z)}{\cdot} f^{i}\colon 
	T^i_{Z/X}\xrightarrow{\quad} \prod_{p}\HH_{(Z)}^{i+p}(X, 
	\Omega_X^{p})\,.
	$$
\end{sit}

\section{Duality and the Infinitesimal Abel-Jacobi Map}
\begin{sit}
	Due to their naturality, the trace maps with support as well as the 
	(powers of the) Atiyah class localize to yield 
	generalized semiregularity maps on the level of sheaves
	$$
			{\tilde \tau}^{i}\colon\cale xt^{i}_X (\calf,\calf)\to 
			\prod_{p}\calh^{i+p}_{Z}\left( X, \Omega^p_X\right)
	$$
	and these morphisms of sheaves fit together to morphisms of the 
	associated local-global spectral sequences
	\begin{diagram}
		\HH^{i}\left(X, \cale xt^{j}_X (\calf,\calf)\right) &\rImplies&
		\Ext^{i+j}_{X}(\calf,\calf)\\
		\dTo^{\dis\HH^{i}\left(X,{\tilde \tau}^{j}\right)}&&
		\dTo^{{\dis\tau}_{Z}^{i+j}}\\
		\dis\prod_{p}\HH^{i}\left(X,\calh^{p+j}_{Z}\left( X, 
		\Omega^p_X\right)\right)
		&\rImplies&
		\dis\prod_{p}\HH_{Z}^{i+j+p}\left(X,\Omega^p_X\right)\,.
	\end{diagram}
As $X$ is assumed to be smooth, each $\Omega^p_X$ is locally free and 
$\calh^{k}_{Z}\left( X, \Omega^p_X\right)=0$ for $k < {\rm 
depth}_{Z}\calox$.  It follows that the components $\tau_{Z}^{i,p}$, 
and then a fortiori $\tau^{i,p}$ and $\sigma_{Z}^{i,p}$ or 
$\sigma^{i,p}$ for $Z$, vanish whenever $i+p < {\rm depth}_{Z}\calox$.
\end{sit}

\begin{sit}
	If $Z\subseteq X$ is a closed subspace of codimension $q$ that is 
	locally Cohen-Macaulay, then the local cohomology sheaves 
	$\calh_{Z}^{q'}(X,\Omega^{p}_X)$ vanish for $q' < q$ and any $p$, 
	and so there are natural inclusions $ 
	\HH^{i}(X,\calh_{Z}^{q}(X,\Omega^{p}_X))\hookrightarrow 
	\HH^{i+q}_{Z}(X,\Omega^{p}_X)$.
	
	Now assume that the given subspace satisfies furthermore $T^2_{Z/X} = 
	\HH^{1}(X,\caln_{Z/X})$.  The isomorphisms of $\calox$-modules
	$$
	\caln_{Z/X} \cong \calh om_{Z}(\cali/\cali^{2},\caloz)\cong \calh 
	om_{X}(\cali,\caloz) \cong \cale xt^{1}_X (\caloz,\caloz)
	$$
	show that the semiregularity maps with support, 
	$\sigma^{2,q-1}_{Z}\colon T^2_{Z/X}\to 
	\HH^{q+1}_{Z}(X,\Omega^{q-1}_X)$, factors through $\HH^{1}(X,{\tilde 
	\sigma}_{Z}^{1,q-1})$.  Thus the semiregularity map $\sigma^{2,q-1}$ 
	itself fits into the following commutative diagram
	\begin{equation*}
		\xymatrix{T^2_{Z/X} =\HH^1(X,\caln_{Z/X})
		\ar[rrr]^{\HH^{1}(X,{\tilde 
		\sigma}_{Z}^{1,q-1})} \ar[d]_{\sigma^{2,q-1}} &&& 
		\HH^{1}\left(X,\calh_{Z}^{q}(X,\Omega^{q-1}_X)\right) \ar[d]^{\rm 
		inclusion}\\
		\HH^{q+1}(X,\Omega^{q-1}_X) &&& \ar[lll]_{\rm canonical} 
		\HH^{q+1}_{Z}(X,\Omega^{q-1}_X)}
	\end{equation*}
\end{sit}

\begin{sit}
	If $Z$ is locally a complete intersection of codimension $q$ in 
	$X$, it is in particular locally Cohen-Macaulay and satisfies 
	$T^2_{Z/X} = \HH^{1}(X,\caln_{Z/X})$.  With $d=\dim Z$ and $\cali$ 
	the defining ideal of $Z$, take exterior powers of the Jacobi map 
	$j:\cali/\cali^2\to \Omega^{1}_X\otimes\caloz$ to define for $0\le 
	i\le q$ natural pairings
	$$
	\Omega_X^{d+i}\otimes \Lambda^{q-i}\left(\cali/\cali^{2}\right)
	\xrightarrow{1\otimes \Lambda^{q-i} j}
	\Omega_X^{d+i}\otimes \Omega_X^{q-i}\otimes \caloz
	\xrightarrow{(\ \wedge\ )\otimes 1}
	\omega_X\otimes\caloz\,.
	$$
	As $\cali/\cali^{2}$ is locally free, these give rise, 
	equivalently, to $\calox$-homomorphisms
	$$
	\gamma_{i}\colon \Omega_X^{d+i}\xrightarrow{\quad}
	\omega_X\otimes\Lambda^{q-i}\caln_{Z/X}
	\cong \omega_Z\otimes\Lambda^{i}\left(\cali/\cali^{2}\right)
	$$
	where we have used the adjunction formula 
	$\omega_Z\cong\omega_X\otimes\det\caln_{Z/X}$ for the dualizing 
	module $\omega_{Z}$, and the isomorphism 
	$\Lambda^{q-i}\caln_{Z/X}\cong 
	\det\caln_{Z/X}\otimes\Lambda^{i}\left(\cali/\cali^{2}\right)$.
	
	The duality theorem for finite morphisms says that
	$$
	\Ext^{p}_{X}(\omega_Z\otimes\Lambda^{i}\left(\cali/\cali^{2}\right),\omega_X)\cong 
	\Ext^{p-q}_{Z}(\Lambda^{i}\left(\cali/\cali^{2}\right),\caloz)\cong
	\HH^{p-q}(Z,\Lambda^{i}\caln_{Z/X})
	$$
	and clearly
	$$
	\Ext^{p}_{X}(\Omega_X^{d+i},\omega_X)\cong
	\HH^{p}(X,\Omega_X^{q-i})
	$$
	for any $p$ or $i$, whence ${\check \gamma}^{i,p} = 
	\Ext^{p}_{X}(\gamma_i,\omega_X)$ can be identified as a map between 
	cohomology groups,
	$$
	{\check \gamma}^{i,p}\colon 
	\HH^{p-q}(Z,\Lambda^{i}\caln_{Z/X})\xrightarrow{\quad}
	\HH^{p}(X,\Omega_X^{q-i})\,.
	$$
	Equivalently, $\check\gamma^{i,p}$ is the Serre dual of the
map $\HH^{d+q-p}_c(X,\Omega_X^{d+i})\to
\HH^{d+q-p}_c(Z,\omega_Z\otimes \Lambda^{i} \cali/\cali^2)$
induced by $\gamma_i$.

\end{sit}
	
\begin{sit}
	Some of these maps are classical objects:
	\begin{enumerate}
		\item For $i=0, p=q$, and with $f:Z\hookrightarrow X$ the closed 
		embedding, one finds by \cite[149-20, Thm.4]{FAG} the leading term 
		of the direct image map $f_{*}$ in cohomology, 
		$$
		{\check \gamma}^{0,q}= 
			f_{*}\left|_{\HH^{0}(Z,\caloz)}\right.\colon 
			\HH^{0}(Z,\caloz)\xrightarrow{\quad} 
			\HH^{q}(X,\Omega_X^{q})\,.
		$$
		By definition, ${\check \gamma}^{0,q}(1_{Z}) = [Z]$, the 
		fundamental (Hodge) class of $Z$, and a simple application of 
		the Grothendieck-Riemann-Roch Theorem, see \cite[2.11]{Murre}, 
		shows ${\check \gamma}^{0,q}(1_{Z}) = \ch_{q}(\caloz)$ which 
		is the result quoted at the beginning (\ref{formula}).
			
		\item For $i=1, p=q$ and $Z$ smooth on a compact K\"ahler 
		manifold, the map ${\check \gamma}^{1,q}$ is the {\em 
		infinitesimal Abel-Jacobi map\/} at $Z$ as shown in 
		\cite[Lect.2, p.29]{Gre},
		$$
		{\check \gamma}^{1,q}=AJ_{X,Z}\colon 
		\HH^0(X,\caln_{Z/X})\xrightarrow{\quad} 
		\HH^{q}(X,\Omega^{q-1}_X)\,.
		$$
			
		\item For $i=1, p=q-1$ and $Z$ again locally a complete 
		intersection, one finds Bloch's semiregularity map
		$$
		{\check \gamma}^{1,q+1}=\sigma_{\rm Bloch}:
		\HH^1(Z,\caln_{Z/X})\xrightarrow{\quad} 
		\HH^{q+1} (X,\Omega^{q-1}_X)\,.
		$$
	\end{enumerate}
\end{sit}
%
%


These observations lead to the following duality theorem for the 
leading terms of the generalized semiregularity maps.

\begin{theorem}
	Let $Z\subseteq X$ be a closed subspace of codimension $q$ that is 
locally a complete intersection.
	\begin{enumerate}
		\item[(1)] One has $T^{p}_{Z/X}\cong 
		\HH^{p-1}(Z,\caln_{Z/X})$ and
		$$
		\sigma^{p,q-1} \cong {\check \gamma}^{1,p+q-1}\colon
		T^{p}_{Z/X}\xrightarrow{\quad} \HH^{p+q-1}(X,\Omega_{X}^{q-1})\,.
		$$
		Each of these maps fits into a commutative diagram
		\begin{equation*}
			\xymatrix{
			\HH^p(X,\caln_{Z/X}) 
			\ar[rrr]^{\HH^{p}(X,{\tilde \sigma}^{1,q-1})} 
			\ar[d]_{\sigma^{p,q-1}}&&& 
			\HH^{p}\left(X,\calh_{Z}^{q}(X,\Omega^{q-1}_X)\right)
			\ar[d]^{\rm inclusion}\\
			\HH^{p+q}(X,\Omega^{q-1}_X) &&& \ar[lll]_{\rm canonical} 
			\HH^{p+q}_{Z}(X,\Omega^{q-1}_X)\,.}
		\end{equation*}
		
		\item[(2)] Bloch's semiregularity map satisfies $\sigma_{\rm 
		Bloch} = \sigma^{2,q-1}$.
	
		\item[(3)] If $Z$ is smooth in a compact K\"ahler manifold $X$, 
		then the infinitesimal Abel-Jacobi map at $Z$ is given by the 
		component $\sigma^{1,q-1}$ of the generalized semiregularity 
		map.		
	\end{enumerate}
\end{theorem}

Clearly $(1)\Longrightarrow (2), (3)$ in view of the explicit 
construction of both $\sigma_{\rm Bloch}$ and $AJ_{X,Z}$.  
Verification of $(1)$ is achieved through a local calculation that 
identifies ${\tilde \sigma}^{1,q-1}$.  The patient reader is referred 
to the details in \cite{BF} whereas the impatient one may wish to 
extract the calculation from the explicit formalism in \cite{ALJ}.\qed

\begin{sit}
	Of course, this result raises several questions, in particular whether for 
	any algebraic cycle $Z$ on a smooth space the diagram
$$
		\xymatrix{ 
		\HH^0(X,\caln_{Z/X})= T^1_{Z/X}
		\ar[d]_{
		\begin{smallmatrix}
			\rm Infinitesimal\\
			\rm Abel-Jacobi\ map
		\end{smallmatrix}
		} 
		\ar[rrr]^{f^{1}}
		&&& 
		\Ext^{1}_{X}(\caloz,\caloz)
		\ar[d]_{\dis\tau_{Z}^{1,q-1}=}^{\dis\tr_{Z}(e^{ -\at_{X}(\caloz)}\smallsmile *)^{q-1}}\\
		\HH^{q}(X,\Omega^{q-1}_X) &&& 
		\ar[lll]_{\rm canonical}
		\HH^{q}_{Z}(X,\Omega^{q-1}_X)}
$$
	commutes. Note that $f^{1}$ is always a split monomorphism by (\ref{Z2}) 
	so that the general infinitesimal Abel-Jacobi map should be a retract of 
	$\tau^{1,q-1}$.
	
	Even more generally, there should be an Abel-Jacobi map for
	deformations of arbitrary coherent sheaves on a compact
	algebraic or K\"ahler manifold such that its differential is
	essentially given by the trace applied to multiplication with powers of the
	Atiyah class. 
\end{sit}

\section{Applications to Deformation Theory}

\begin{sit}
	Ideally, one would like to interpret the family of semiregularity maps 
	$\sigma$ or $\tau$ as maps between cohomology groups of cotangent 
	complexes induced by morphisms of deformation theories.  Whereas the 
	domain of these maps is of the desired form, it is not clear what kind 
	of deformation theory, if any, would allow the groups 
	$\bbbh^{q}:=\prod_{p}\HH^{p+q}(X,\Omega^{p}_{X})$ as its tangent 
	cohomology.  Hodge theory shows that if such a deformation theory 
	existed, it were necessarily non obstructed for any complex manifold 
	$X$ that is bimeromorphically equivalent to a K\"ahler manifold: 
	According to Deligne \cite[(5.5)]{Del}, each of the functors
	$$
	\caln\mapsto \HH^{q}(X_{T},\Omega^{p}_{X_{T}/T}\otimes_{\calo_{T}}\caln)
	$$
	on the coherent modules $\caln$ over an artinian complex germ $T$ 
	is then {\em exact\/}, and this applied to $\bbbh^{1}$ would show 
	that every ``infinitesimal deformation'' could be lifted.  This 
	line of reasoning applies indeed directly in some special cases, see 
	\cite{Ran4}.
	
	In general, even though there might not be an actual deformation 
	theory underlying the groups $\bbbh^{q}$, it was first pointed out 
	by Ran, \cite{Ran1, Ran2, Ran3}, and then made more precise in 
	\cite{Kaw1,Kaw2, Kaw3} and \cite{FM}, that Deligne's theorem can 
	be used to obtain results on embedded deformations or deformations 
	of coherent modules.  Complete proofs and further generalizations 
	of the following results will be contained in \cite{BF}.
\end{sit}

\begin{theorem}
Assume that $f:X\to \Sigma$ is a smooth morphism of complex spaces 
with $X_0:=f^{-1}(0)$ a compact manifold bimeromorphically equivalent 
to a K\"ahler manifold.  Let $\calf_0$ be a coherent module on $X_0$ 
and 
$$
\tau^{2} = \tr(e^{-\at(\calf_0)}\smallsmile *)\colon 
\Ext^2_{X_0} (\calf_0, \calf_0) \xrightarrow{\quad}
\prod_{p \ge 0} H^{p+2}
(X_0, \Omega^p_{X_0})
$$
the corresponding generalized semiregularity map.
The dimension of the basis $S$ of a semi-universal deformation of $\calf_0$
satisfies
 $$\dim S \ge \dim_\bbbc \Ext^1_{X_0}(\calf_0, \calf_0) -
\dim_\bbbc \ker \tau^{2}\,.$$
If $\tau^{2}$ is injective then $S$ is smooth over a closed subspace 
of $\Sigma$.  \qed
\end{theorem}

\begin{sit}
	This result generalizes in characteristic zero the {\em 
	Artamkin-Mukai criterion\/} \cite{Arta1, Mukai, Kaw1} which is the 
	above statement but restricted to the component 
	$\tau^{2,0}:\Ext^2_{X} (\calf, \calf) \xrightarrow{\quad} H^{2}(X, 
	\calo_{X})$ in the absolute case where $X=X_{0}$ is a nonsingular 
	projective variety.  The map $\tau^{2,0}$ is just the ordinary 
	trace and this component of the generalized semiregularity map has 
	indeed a description as a transformation between tangent 
	cohomologies induced by a morphism of deformation theories: As $X$ 
	is smooth, every coherent sheaf $\calf$ admits a well defined 
	determinant line bundle $\det\calf$, and associating to $\calf$ 
	its determinant is a morphism of deformation functors that induces
	$\tau^{i,0}\colon \Ext^i_{X} (\calf, \calf) \xrightarrow{\quad} H^{i}(X, 
	\calo_{X})$ on the corresponding tangent cohomology groups.
	
	As the target, the Picard functor of $X$, is smooth, it follows 
	that all obstructions to deform $\calf$ are already contained in 
	the kernel of the trace map $\tau^{2,0}$.  This deformation 
	theoretic interpretation allows Artamkin to obtain as well results 
	in positive characteristic that are beyond reach here.
\end{sit}

The analogous application to embedded deformations of subspaces generalizes 
Bloch's original result \cite[Thm.(7.3)]{Bloch} relating 
semiregularity and smoothness of the Hilbert scheme.  It encompasses 
the previous generalizations of that result in \cite{Ran2, Kaw1}.

\begin{theorem}
Assume that $f:X\to \Sigma$ is a smooth morphism of complex spaces and 
that $X_0:=f^{-1}(0)$ is a compact manifold bimeromorphically 
equivalent to a K\"ahler manifold.  Let $Z\subseteq X_0$ be a closed 
subspace and
$$
\sigma^{2}\colon
T^2_{Z/X_0} (\caloz) \xrightarrow{\quad} \prod_{p \ge 0} H^{p+2}
(X_0, \Omega^p_{X_0})\ .
$$
the generalized semiregularity map.

The dimension of the Douady space $H_{X/\Sigma}$ at $[Z]$ satisfies
$$
\dim_{[Z]} H_{X/\Sigma} \ge \dim_\bbbc T^1_{Z/X_0}(\caloz) -
\dim_\bbbc \ker \sigma^{2}\,.
$$

If $\sigma^{2}$ is injective, then $H_{X/\Sigma}$ is
smooth over a closed subspace of $\Sigma$ in a neighbourhood of
$[Z]$.
\qed
\end{theorem}

\section{The Singular Case}
In this announcement of results we restricted ourselves to smooth 
spaces $X$, but all the necessary ingredients
\begin{enumerate}
	\item[$\bullet$] Description of the obstruction for embedded 
	deformations or for deformations of complexes of modules in terms 
	of the Atiyah class,

	\item[$\bullet$] Construction of a trace for perfect complexes,

	\item[$\bullet$] Construction of the Atiyah-Chern character for perfect complexes,
\end{enumerate}
have been established by L.~Illusie in arbitrary generality in the 
algebraic case, and the analytic case can be handled through the 
methods of \cite{Fle, OTT}. 

Using this machinery, one constructs first the cotangent complex 
$\bbbl_{X/\Sigma}$ for an arbitrary morphism $X\to \Sigma$ of analytic 
spaces, schemes or even ringed topoi, secondly one defines traces 
$\Ext^{q}(\calf,\calf\otimes^{\bbbl}\Lambda^{p}\bbbl_{X/\Sigma})\to 
\prod_{p}\HH^{p+q}(X, \Lambda^{p}\bbbl_{X/\Sigma})$ for a perfect 
complex $\calf$, which needs special care as 
$\Lambda^{p}\bbbl_{X/\Sigma}$, the derived exterior powers of the 
cotangent complex, will in general not be perfect, and finally one 
defines the Atiyah class and, in characteristic zero, an 
Atiyah-Chern-Illusie character in 
$\prod_{p}\Ext^{p}(\calf,\calf\otimes^{\bbbl}\Lambda^{p}\bbbl_{X/\Sigma})$ 
to obtain a commutative diagram extending (\ref{key}) to something 
resembling a metatheorem:
\begin{theorem}
	Let $X\to \Sigma$ be a morphism of complex spaces or of schemes over 
	a field of characteristic zero. If $\calf$ is a
	perfect complex on $X$ and $\caln$ is a coherent
	$\calox$-module then there is a commutative diagram
	$$
	\xymatrix{
	T^1_{X/\Sigma}(\caln) = \Ext^{1}_{X}(\bbbl_{X/\Sigma},\caln)
	\ar[rrr]^{\quad {\rm ob}_{\calf}=*\lr -\at_{X/\Sigma}(\calf)}
	\ar[dd]_{*\lr \ch_{X/\Sigma}(\calf)} &&&
	\Ext^2_{X}(\calf, \calf\otimes^{\bbbl} \caln)
	\ar[llldd]^{\qquad\tr(e^{-\at_{X/\Sigma}(\calf)}
	\otimes ^{\bbbl} \caln\smallsmile *)}\\
	\\
	\dis\prod_{p} H^{p+2}(X,\Lambda^p(\bbbl_{X/Y})\otimes^{\bbbl}\caln)&.}
	$$
	Moreover, if $\calf=\caloz$ for $Z\subseteq X$ a closed subspace,
	then the morphism on top, that is, contracting against the
	opposite of the Atiyah class, factors naturally through
	$T^{2}_{Z/X}(\caln|_{Z})=\Ext^{2}_{Z}(\bbbl_{Z/X},\caloz\otimes^{\bbbl}\caln)$.
	\qed
\end{theorem}


\end{document}